\documentclass{amsart}

\usepackage{amssymb}

\theoremstyle{plain}

\numberwithin{equation}{section}

\newtheorem*{teo}{Theorem A}
\newtheorem*{teoremab}{Corollary}

\newtheorem{cora}[equation]{Corollary}

\newtheorem{lem}[equation]{Lemma}

\newtheorem{prop}[equation]{Proposition}

\theoremstyle{definition}

\begin{document}	
\title{Homogeneous products of conjugacy classes}

\author{Edith Adan-Bante}

\address{University of Southern Mississippi Gulf Coast, 730 East Beach Boulevard,
 Long Beach MS 39560}

\email{Edith.Bante@usm.edu}


\subjclass{20d15}

\date{2005}
\begin{abstract} Let $G$ be a finite group and $a\in G$. Let 
$a^G=\{g^{-1}ag\mid g\in G\}$ be the conjugacy class of $a$ in $G$. 
Assume that $a^G$ and $b^G$ are conjugacy classes of $G$ 
with the property that ${\bf C}_G(a)={\bf C}_G(b)$. Then $a^G b^G$ is 
a conjugacy class if and only if $[a,G]=[b,G]=[ab,G]$ and 
$[ab,G]$ is a normal 
subgroup of $G$. 
\end{abstract}
\maketitle

\begin{section}{Introduction}

Let $G$ be a finite group, $a\in G$  and $a^G=\{a^g\mid g\in G\}$ 
be the conjugacy class of $a$ in $G$. Denote by $|a^G|$ the size of the 
set $a^G$. 
Given $a,g\in G$, set $[a,g]=a^{-1}a^g$. Also set $[a,G]=\{[a,g]\mid g\in G\}$. 
Let ${\bf C}_G(a)=\{ g\in G\mid a^g=a\}$ be the centralizer
of $a$ in $G$ and $1_G$ be the identity of $G$.
 Through this note, we will use the well known fact that 
$|a^G|=|G:{\bf C}_G(a)|$.

In Theorem A of  \cite{edithma}, 
it is proved that if $G$ is a finite nilpotent group and 
$\chi, \psi$ are faithful irreducible characters with the property that
$\chi\psi$ is a multiple of an irreducible, then $\chi$ and $\psi$ both  vanish outside
the center ${\bf Z}(G)$ of $G$, i.e. $\chi(g)=\psi(g)=0$ for all $g\in G\setminus {\bf Z}(G)$. This note was motivated 
by wondering what would be the 
analogous result in conjugacy 
classes.

Let $a^G$ and $b^G$ be conjugacy classes such that the product
$a^G b^G= \{xy\mid x\in a^G, y\in b^G\}$  is also a conjugacy class.
We can check that $(ab)^G$ is a subset of $a^G b^G$ and thus if
 $a^G b^G$ is a conjugacy class, then $a^G b^G=(ab)^G$. 
Is there any relationship between $a$, $b$ and $G$? The answer
in general seems to be no. For instance, if we take   any group $G$,
any element $a$ of $G$, then $a^G 1_G^G= (a1_G)^G$.  
But if we add the additional hypothesis that 
 ${\bf C}_G(a)={\bf C}_G(b)$, then we have the following
 
\begin{teo}\label{eta1prop}
Let $G$ be a finite group, $a^G$ and $b^G$ be  conjugacy classes of $G$. Assume that 
${\bf C}_G(a)={\bf C}_G(b)$.
 Then  $a^G b^G=(ab)^G$ if and only if  $[ab,G]=[a,G]=[b,G]$ and $[ab,G]$ is a normal
subgroup of $G$. 
In particular, given any conjugacy class $a^G$ of $G$, then 
 $a^G a^G=(a^2)^G$ if and only if $[a,G]$ is a normal subgroup of $G$.
 \end{teo}
 
%Let $S_3$ be the symmetric group of 3 letters and $a$ any element of order 3.
%e can check that $[a,G]$ is not a subgroup.
%ecause for any finite group $G$ and any $a$ any element of $G$, 
%e have that $a^G 1^G= (a1_G)^G$, the hypothesis that ${\bf C}_G(a)={\bf C}_G(b)$ is 
%quired in order to obtain that $[ab,G]$ is a subgroup.

  We regard the hypothesis that $a^G$ and $b^G$ are conjugacy classes of $G$ with 
${\bf C}_G(a)={\bf C}_G(b)$ as the ``dual" to the hypothesis that 
two irreducible characters have the same kernel.   

  The following  is a direct application of Theorem A.
 
 \begin{teoremab}
 Let $G$ be a finite nonabelian simple group, $a^G$ and 
 $b^G$ be  conjugacy classes of $G$. Assume that 
${\bf C}_G(a)={\bf C}_G(b)$. Then $a^G b^G=(ab)^G$ if and only if $a=b=1_G$.
In particular, $a^G a^G=(a^2)^G$ if and only if  $a=1_G$.
 \end{teoremab}

   Is it possible to find a finite group $G$ and a conjugacy class $a^G$ of $G$ such
   that $a^G a^G=(a^2)^G$ and $|a^G|=2$? The answer is no, such group with such 
   conjugacy class can not exist. In 
    Proposition \ref{squaresofsize2}, we show 
    that if $G$ is a finite group, $a^G$ and $b^G$ are
    conjugacy classes such that  ${\bf C}_G(a)={\bf C}_G(b)$, 
    then necessarily
    $a^G b^G$ is the union of exactly 2 distinct conjugacy classes.
    But then, is it possible to  find a finite group $G$ and 
    a conjugacy class $a^G$ of $G$ such
   that $a^G a^G=(a^2)^G$ and $|a^G|$ is a power of 2? If, in addition, we require 
   that
   the group is supersolvable, then the answer is again no.  More specifically, 
   in Proposition \ref{supersolvable} is shown that
    if $G$ is a finite supersolvable
  group, $a^G$ and $b^G$ are conjugacy classes of $G$ with 
  ${\bf C}_G(a)={\bf C}_G(b)$ and $|a^G|=2^n$ for some integer $n>0$, then $a^G b^G$ is 
  the union of at least 2 distinct conjugacy classes. The author wonders if the answer
   remains no if we do not require that the group is supersolvable. 
  On the other hand, 
  in Proposition \ref{eta1example},  given any odd integer
  $n>0$,
 we provide an example of
 a nilpotent group $G$ and a conjugacy class $a^G$ such that 
 $a^G a^G= (a^2)^G$  and $|a^G|=n$. 

We want to close this introduction by mentioning that there is a number of papers
concerning   
 products of conjugacy classes and finite 
 groups. A very recent development  is  \cite{dade}, where the authors classify all finite
 groups $G$ such that the product of any two non-inverse conjugacy classes of $G$ is 
 always a conjugacy class of $G$.

{\bf Acknowledgment.} I would like to thank Professor Everett C. Dade for his suggestions
 to improve both the results and the presentation of this note. I also thank the referee
 and the editor for useful comments.   
\end{section}

\begin{section}{Proof of Theorem A}

 We will denote by $1_G$ the identity of the group $G$. 
\begin{lem}\label{productsn} Let $G$ be a finite group and  $a,b\in G$. Then
\begin{equation*}\label{asquare}
a^G b^G =ab[a^{b^{-1}},G][b,G].
\end{equation*}
Thus if $a^b=a$ then $a^G b^G=ab[a,G][b,G]$.
\end{lem}
\begin{proof}
Observe that  
\begin{equation*}
\begin{split}
a[a,G]b[b,G]& = ab[a,G]^{b^{-1}}[b,G] \\
&= ab[ a^{b^{-1}}, G^{b^{-1}}][b,G]\\
& = ab [a^{b^{-1}},G] [b,G].
\end{split}
\end{equation*}
\end{proof}
\begin{lem}\label{groupimpliesnormal}
Let $G$ be a finite group and $c\in G$. If $[c,G]$ is a subgroup of $G$, then 
$[c,G]$ is a normal subgroup of $G$. 
\end{lem}
\begin{proof}
Let $g\in G$ and $x\in [c,G]$. By definition, $(c)^g=cy$ for some $y \in [c,G]$, and 
$cx=(c)^h$ for some
$h\in G$. Also $(c)^{hg}= cw$ for some $w\in [c,G]$.
Observe that  $(c)^{hg}= ((c)^h)^g= (cx)^g= (c)^gx^g= cyx^g$. Thus $yx^g=w$ and 
$x^g=y^{-1}w\in [c,G]$.
We conclude that $[c,G]$ is a normal subgroup of $G$. 
\end{proof}
 \begin{proof}[Proof of Theorem A]
 Since ${\bf C}_G(a)={\bf C}_G(b)$, we have that $ab=ba$.
 Observe that if $[a,G]=[b,G]=[ab,G]$ and $[ab,G]$ is a normal subgroup, then by Lemma 
 \ref{productsn}, we have that $a^G b^G= ab[a,G][b,G]=ab[ab,G]=(ab)^G$. 
 We may assume
 now that $a^Gb^G=(ab)^G$ and we want to conclude that $[a,G]=[b,G]=[ab,G]$ and $[ab,G]$ is a normal subgroup of $G$. 
 
Since ${\bf C}_G(a)={\bf C}_G(b)$, we have that $|a^G|=|G:{\bf C}_G(a)|=|G:{\bf C}_G(b)|=|b^G|$,
  ${\bf C}_G(ab)\geq {\bf C}_G(a)$ and therefore
$|(ab)^G|\leq |a^G|$. Because $a^G b^G=(ab)^G$, we have then that 
$|(ab)^G|=|a^G|=|b^G|$. Since $a^G b^G=(ab)^G$ and $ab=ba$,
by Lemma \ref{productsn} we have that 
$[a,G][b,G]=[ab,G]$. 
Thus $[a,G]=[b,G]=[ab,G]$ since $|(ab)^G|=|a^G|=|b^G|$ and $1_G$ is in both $[a,G]$ and $[b,G]$.

Since $[a,G]=[b,G]=[ab,G]$ and $[a,G][b,G]=[ab,G]$, then $[ab,G][ab,G]=[ab,G]$. Clearly
$[ab,G]$ is nonempty since 
$1_G\in [ab,G]$. 
 We conclude that $[ab,G]$ is a subgroup of $G$ since   
 $uv\in [ab,G]$ for any $u,v$ in $[ab,G]$ and $[ab,G]$ is a nonempty
 finite set. The result then follows by Lemma \ref{groupimpliesnormal}.   
\end{proof}
\end{section}

\begin{section}{Further Results}  
%{\bf Notation}. 
%Let $G$ be a finite group and $N$ be a normal subgroup of $G$. Let $\gamma_N: G\mapsto G/N$
%be the epimorphism $\gamma_N(a)= aN$.

Let $X$ be a $G$-invariant subset of $G$, i.e. 
$X^g=\{x^g\mid x\in X\}=X$ for all $g\in G$.
  Then $X$ can be expressed as a union of 
  $n$ distinct conjugacy classes of $G$, for some integer $n>0$. Set
 $\eta(X)=n$.
 
\begin{lem}\label{lemma1}
Let $G$ be a finite $p$-group
 and $N$ be a normal subgroup of $G$. 
Let $a$ and $b$ be elements
of $G$. 
If $(aN)^{G/N}\cap (bN)^{G/N}= \emptyset$
 then $a^G\cap b^G= \emptyset$.
 Thus $\eta((aN)^{G/N} (bN)^{G/N})\leq 
  \eta(a^G b^G )$.
\end{lem}
\begin{proof}
See Lemma 2.1 of  \cite{edith2}.
\end{proof}

\begin{prop}\label{notinthecenter}
Let $G$ be a group of odd order and $a^G$ be the conjugacy class of $a$ in $G$.
Then 
\begin{equation}\label{centerg}
{\bf Z}(G)\cap a^G a^G\neq \emptyset
\end{equation}
\noindent if and only if $|a^G|=1$.
Thus if $|a^G|>1$ and  $b^G\subseteq a^Ga^G$, then $|b^G|>1$.
\end{prop}
\begin{proof}
Suppose that there exist some $z\in {\bf Z}(G)\cap a^G a^G$. Then
there exist some $g\in G$ such that $aa^g=z$. Thus $a^g=a^{-1} z$
and therefore $(a^{-1})^g=a z^{-1}$. 
Observe that 
\begin{equation*}
a^{g^2}= (a^g)^g= (a^{-1} z)^g= (a^{-1})^g z= (a z^{-1})z=a.
\end{equation*}
Thus $g^2\in {\bf C}_G(a)$. Since $G$ is of odd order, 
$g^2\in {\bf C}_G(a)$ implies that   $g\in {\bf C}_G(a)$.
So  $a^2= z$ and $a \in {\bf Z}(G)$. We conclude that $|a^G|=1$. 
\end{proof}

Let $E$ be an extraspecial group of order $3^3$ and exponent $3$. Let $a\in E\setminus {\bf Z}(E)$. Set $b=a^2$. We can check that ${\bf C}_G(a)={\bf C}_G(b)$
and $a^G b^G = {\bf Z}(E)$. Thus given a finite group $G$ of odd order, conjugacy classes 
$a^G$ and $b^G$ of $G$ with ${\bf C}_G(a)={\bf C}_G(b)$,   
${\bf Z}(G)\cap a^G b^G\neq \emptyset$ may not imply that $|a^G|=1$. 

Let $Q_8$ be the quaternion group and $a\in Q_8$ be an element of order 4. 
We can check that $a^{Q_8}=\{a, a^{-1}\}$ and $a^{Q_8}a^{Q_8}={\bf Z}(Q_8)$.
Thus Proposition \ref{notinthecenter} may not remain true if the 
group $G$ has even order.

 \begin{prop}\label{squaresofsize2}
 Let $G$ be a finite group, $a^G$ and $b^G$ be conjugacy classes  with 
 ${\bf C}_G(a)={\bf C}_G(b)$
 and
 $|a^G|=2$.
Then $\eta(a^G b^G)= 2$.
In particular 
$\eta(a^Ga^G)=2$.
\end{prop}
\begin{proof}
Set $N={\bf C}_G(a)$. Observe that $N$ is a normal subgroup of $G$ since 
$|G:N|=|G:{\bf C}_G(a)|=2$.
Since ${\bf C}_G(a)={\bf C}_G(b)$, we have that $ab=ba$. Fix $g\in G\setminus N$.
Since $|a^G|=|b^G|=2$, ${\bf C}_G(a)={\bf C}_G(b)$ and $g\in G\setminus N$, 
we have that
 $a^G=\{a,a[a,g]\}$ and $b^G=\{b,b[b,g]\}$. Also $[a,g]b=b[a,g]$ since $[a,g]\in N={\bf C}_G(b)$.
Therefore
\begin{equation}\label{productab}
a^G b^G=\{ab, ab[a,g], ab[b,g], ab[a,g][b,g]\}. 
\end{equation}

Since $|G:N|=2$, we have that
$a^{g^2}=a$ and $a^{g^2}=(a[a,g])^g= a^g[a,g]^g=a[a,g][a,g]^g$. Thus $[a,g][a,g]^g=1_G$. 
  Fix  $n\in N$. Observe that 
\begin{equation}\label{first}
a^{gn}= (a^g)^n=(a[a,g])^n=a^n [a,g]^n=a[a,g]^n.
\end{equation}
Also observe that
\begin{equation}\label{2}
a^{ng[g,n]} = (a^n)^{g[g,n]}= a^{g[g,n]} = (a[a,g])^{[g,n]}= a^{[g,n]} [a,g]^{[g,n]}=a[a,g]^{[g,n]}.
\end{equation}
Since $gn=ng[g,n]$, we have that $a^{gn}=a^{ng[g,n]}$. Thus 
 by \eqref{first} and \eqref{2} we have that $[a,g]^n=[a,g]^{[g,n]}$
and so $[a,g]= [a,g]^{(n^{-1})^g}$. Thus $N^g\in {\bf C}_G([a,g])$. 
 Since $N$ is normal in $G$, we conclude that 
$[a,g]^n=[a,g]$ for any $n\in N$. Similarly, we can check that 
$[b,g]^n=[b,g]$ for any  $n\in N$.

 Since $[a,g][a,g]^g=1_G$, $[a,g]^n=[a,g]$ and $[b,g]^n=[b,g]$ for all $n\in N$, and
 $ab=ba$, we have
 \begin{equation*}
 \begin{split} 
  (ab[a,g])^g & = (ba[a,g])^g=b^g a^g [a,g]^g\\
  &= b[b,g]a[a,g][a,g]^g= ba[b,g][a,g][a,g]^g \\
  &=ab [b,g].
  \end{split}
  \end{equation*}
 Thus $(ab[a,g])^G=\{ab[a,g], ab[b,g]\}$ since $[a,g]^n=[a,g]$ and $[b,g]^n=[b,g]$
  for all $n\in N$, $|G:N|=2$ and 
 $g\in G\setminus N$.
 
 Since $(ab)^G=\{ab, ab[a,g][b,g]\}$, $[a,g]\neq 1$ and $[b,g]\neq 1$, we conclude 
 that $\{ab,ab[a,g]b,g]\}$ and
 $ \{ab[a,g],ab[b,g]\} $ are two distinct conjugacy classes. By \eqref{productab}
 we have then that  $a^G b^G=\{ab,ab[a,g][b,g]\} \cup \{ab[a,g],ab[b,g]\} $.
  Therefore
 $\eta(a^G b^G)=2$. 
 \end{proof}
  
\begin{prop}\label{supersolvable}
Let $G$ be a finite supersolvable group, $a^G$ and $b^G$ be conjugacy classes of $G$ with 
${\bf C}_G(a)={\bf C}_G(b)$ and $|a^G|=2^n$ for some integer $n>0$.
Then $\eta(a^G b^G)\geq 2$. 
\end{prop}
\begin{proof}
Let $G$ be a supersolvable group, $a^G$ and $b^G$ be conjugacy classes with 
${\bf C}_G(a)={\bf C}_G(b)$ and $|a^G|=2^n$ for some integer $n>0$.
Assume that the statement is false and  $G$, $a^G$ and $b^G$ are a 
minimal counterexample of the statement
with respect to  $|a^G|$, i.e.  $\eta(a^G b^G)=1$ and for any supersolvable 
group $K$, any 
conjugacy classes $c^K$ and $d^K$ of $K$ such 
that ${\bf C}_K(c)={\bf C}_K(d)$ and $|c^K|=2^m$, where $0<m<n$, then necessarily
 $\eta(c^K d^K)\geq 2$. We are going to conclude that then $2\leq 1$ and thus the statement 
 holds.
 
 By Theorem A, the set $[ab,G]$ is a normal subgroup 
of $G$. Let $N$ be a normal subgroup of $G$ contained in $[ab,G]$ with $|N|=2$. 
Observe such 
subgroup exists since $[ab,G]$ is a normal subgroup of $G$ and $G$ is supersolvable.
Consider now the group $G/N$. Since $N\leq [ab,G]$ and $|N|=2$, we can check that 
$|(aN)^{G/N}|=|a^G|/2$.  Since  $\eta(a^G b^G)=1$, 
by Proposition \ref{squaresofsize2} we have that $|a^G|=2^n>2$, that is $n>1$,
and therefore $|(aN)^{G/N}|=2^{n-1}>1$, i.e. $n-1>0$. Observe that
${\bf C}_{G/N} (aN)={\bf C}_{G/N} (bN)$ because ${\bf C}_G(a)={\bf C}_G(b)$.
Since $|(aN)^{G/N}|=2^{n-1}<|a^G|$ with $n-1>0$,
 we have that $\eta((aN)^{G/N} (bN)^{G/N})\geq 2$.
By Lemma \ref{lemma1} we have that 
$\eta((aN)^{G/N} (bN)^{G/N})\leq \eta(a^G b^G)=1$ and thus $2\leq 1$. 
\end{proof}

\begin{cora}\label{nilpotent}
Let $G$ be a finite nilpotent group and $a^G$ be a conjugacy class $G$. If  
 $a^G a^G= (a^2)^G$, then $|a^G|$ is an odd number. 
\end{cora}
\end{section}
\begin{section}{Examples}
\begin{lem}\label{observation1}
Let $G$ and $K$ be finite groups, $a^G$ be the conjugacy class of $a$ in $G$ and 
$b^K$ be the conjugacy class of $b$ in $K$. Assume that $\eta(a^G a^G)=1$ and 
$\eta(b^K b^K)=1$. Let $G\times K$ be the direct product of $G$ and $K$. Then 
$\eta((a,b)^{G\times K} (a,b)^{G\times K})=1$ and $|(a,b)^{G\times K}|=|a^G| |b^K|$. 
\end{lem}
\begin{proof} By definition of direct product, we have that
$(a,b)^{G\times K}=\{(x,y)\mid x\in a^G, y\in b^K\}$ and 
$(a,b)^{G\times K} = \{(xu, yv)\mid x,u \in a^G \mbox{ and }y, v\in b^K\}$.
Thus $|(a,b)^{G\times K}|=|a^G| |b^K|$ and $\eta((a,b)^{G\times K} (a,b)^{G\times K})=1$.
\end{proof}
\begin{lem}\label{observation2}
Let $p$ be a prime number and $n>0$.
Let  $E$ be an extraspecial group of order $p^3$ and exponent $p$.
Let $G=E\times E\cdots\times E$ be the direct product of $n$ copies of $E$.
Fix $a=(e_1,e_2,\ldots,e_p)\in G$, where $e_i\in E\setminus {\bf Z}(E)$
for $i=1,\ldots p$.
Then 
\begin{equation}\label{center}
a^G= \{az\mid z\in {\bf Z}(G)\}.
\end{equation}
\noindent Thus $|a^G|=p^n$ 
and $a^G a^G=(a^2)^G$.
Therefore given any prime $p$ and any integer $n> 0$, there exist a $p$-group
$G$ and a conjugacy class $a^G$ such that $a^G a^G=(a^2)^G$ and
$|a^G|=p^l$.
\end{lem}
\begin{proof}
Since $G$ is the direct product of $n$ copies of $E$, then 
${\bf Z}(G)={\bf Z}(E)\times\cdots\times {\bf Z}(E)$ and thus 
$|{\bf Z}(G)|=p^n$.
We can check that given any $z\in {\bf Z}(G)$, there exist some 
$g\in G$ such that $a^g=az$. Also, given any $g\in G$, there exists some
$z\in {\bf Z}(G)$ such that $a^g= az$. Thus \eqref{center} follows and 
the proof is now complete.
\end{proof}
 
\begin{prop}\label{eta1example} 
Given  any odd integer $n\geq 1$, there exist a nilpotent group
$G$ and a conjugacy class $a^G$ such that $a^G a^G=(a^2)^G$ and
$|a^G|=n$.
\end{prop}
\begin{proof}
It follows from Lemmas \ref{observation1} and \ref{observation2}. 
\end{proof}

 \end{section}

\end{document}